\documentclass[12 pt]{article}

\usepackage{pstricks}
\usepackage{pstricks-add}
\usepackage{pst-plot}

\usepackage{amsfonts}
\usepackage{amsmath}
\usepackage{amssymb}
\usepackage{graphicx}
\usepackage{aeguill}

\newtheorem{theorem}{Theorem}
\newtheorem{proposition}{Proposition}
\newtheorem{corollary}{Corollary}
\newtheorem{lemma}{Lemma}
\newtheorem{definition}{Definition}

\def \ep {\hspace{0.2 cm}  $\blacksquare$}
\begin{document}

\title{On products of skew rotations.}
\author{M.D. Arnold  \thanks{Institute for Information Transmission Problems of the Russian Academy of Sciences (Kharkevich Institute), Bolshoi Karetny per. 19, Moscow, 127994, Russia} \thanks{International Institute of Earthquake Prediction Theory and Mathematical Geophysics of the Russian Academy of Sciences, Profsoyuznaya str., 84/32, Moscow, 117997, Russia}\and E.I. Dinaburg \footnotemark[1] \thanks{Schmidt Institute of Physics of the Earth of the Russian Academy of Sciences, B. Gruzinskaya str., 10,  Moscow, 123995, Russia}\and G.B. Dobrushina  \footnotemark[1]\and S.A. Pirogov \footnotemark[1]  \and A.N. Rybko \footnotemark[1]}
\maketitle

\begin{abstract}

Let $\{S_1^t\},\ldots,\,\{S_n^t\}$ be the one-parametric groups of shifts along the orbits of Hamiltonian systems generated by time-independent Hamiltonians $H_1,\ldots,\, H_n$ with one degree of freedom. In some problems of population genetics there appear planar transformations having the form $S^{h_n}_n\cdots S_1^{h_1}$ under some conditions on Hamiltonians $H_1,\ldots,\,H_n$. In this paper we study asymptotical properties of trajectories of such transformations. We show that under classical non-degeneracy condition on the Hamiltonians the trajectories stay in the invariant annuli for generic combinations of lengths $h_1$,..., $h_n$, while for the special case $h_1+\dots+h_n=0$ there exists a trajectory escaping to infinity.  

\end{abstract}

\small{\noindent Mathematics Subject Classification: 37J40, 37J15, 37M05;}
\section{Introduction.}

Let $H_1(p,q)$, $H_2(p,q)$ be two time-independent Hamiltonians with one degree of freedom.
Typical examples can be $H_1=\sqrt{(p-p_1)^2+(q-q_1)^2}$ and $H_2=\sqrt{(p-p_2)^2+(q-q_2)^2}$ where $p_1$, $q_1$, $p_2$, $q_2$ are some constants. 
Sometimes instead of symplectic coordinates $(p,q)$ we shall write $(x,y)$ or in complex notation $z=x+iy$.
 Denote by $\{S_1^t\}$, $\{S_2^t\}$ the one-parametric groups of shifts along the trajectories of the first and second system. Ya.G. Sinai formulated a general question about asymptotic properties of transformations  $T^{(h_1,h_2)}=S_2^{h_2}\cdot S_1^{h_1}$, where $h_1$, $h_2$ are fixed numbers. 

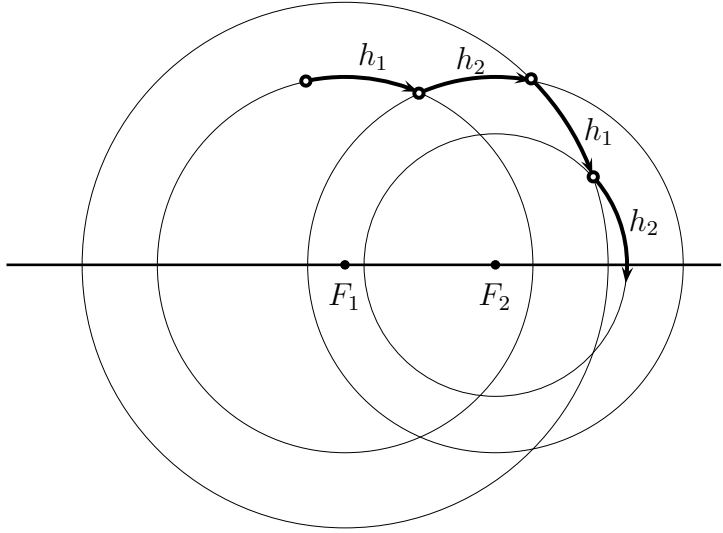
\begin{figure}[ht]
\label{fig: Rotations}
\begin{center}
\begin{pspicture}(0,0)(7,7)
\psline[linewidth=1 pt,linecolor=black,linestyle=solid](-2,3.5)(7.5,3.5)

\psdot(2.5,3.5)
\psdot(4.5,3.5)
\put(2.5,3.5){\pscircle[linewidth= 0.2 pt,linestyle=solid]{2.5}}
\put(2.5,3.5){\pscircle[linewidth= 0.2 pt,linestyle=solid]{3.5}}
\put(4.5,3.5){\pscircle[linewidth= 0.2 pt, linecolor=black,linestyle=solid]{2.5}}
\put(4.5,3.5){\pscircle[linewidth= 0.2 pt,linestyle=solid]{1.75}}

\psarc[linewidth= 1.5 pt, linecolor=black]{<-o}(2.5,3.5){2.5}{67}{102}
\uput[12](2.5,6.2){$h_1$}

\psarc[linewidth= 1.5 pt, linecolor=black]{<-o}(4.5,3.5){2.5}{79}{114}
\uput[12](3.8,6.15){$h_2$}

\psarc[linewidth= 1.5 pt, linecolor=black]{<-o}(2.5,3.5){3.5}{20}{45}
\uput[12](5.5,5.2){$h_1$}

\psarc[linewidth= 1.5 pt, linecolor=black]{<-o}(4.5,3.5){1.75}{-8}{42}
\uput[12](6.1,4){$h_2$}


\uput[12](2.1,3){$F_1$}
\uput[12](4.1,3){$F_2$}


\end{pspicture}
\end{center}
\caption{Product of two skew rotations around the points $F_1$ and $F_2$ with $h_1,\,h_2>0$.}
\end{figure}

Similar transformations appear in some problems of population genetics (see \cite{Sato}).  In the example described above, the map $T^{(h_1,h_2)}$ moves each point $(p,q)$ along the circle with center $(p_1,q_1)$ for the distance $h_1$ in the positive direction and then moves its image along the circle centred at $(p_2,q_2)$ for the distance $h_2$ (see Fig. \ref{fig: Rotations}). For negative values of $h_i$, we assume the movement in the opposite direction.

One can ask a natural question about boundedness of the 
trajectories $\left(T^{(h_1,h_2)}\right)^n(p,q)$.
In more detail: do invariant curves of the map $T^{(h_1,h_2)}$ exist and separate the phase space on bounded annuli?  From our results it follows that for $h_1+h_2\ne 0$ the answers to these questions are affirmative, while for $h_1+h_2=0$ the answers are negative.

%

 
 Assume that in the previous example $p_1=q_1=0$. Then $H_1=\sqrt{p^2+q^2}$ and the corresponding Hamiltonian flow has the form  
  \begin{equation}
\label{eq: simple_shift}
\begin{cases}
\varphi_t=\varphi_0+ \dfrac t \rho\\
\rho_t=\rho_0
\end{cases}
\end{equation}
 in the  usual polar coordinates $\varphi=\mathrm{arg}(z)$, $\rho=|z|$.
For any $t$, this map is symplectic, i.e., it preserves the area element $dx\wedge dy=\rho d\rho\wedge d\varphi$.

\begin{definition}
For a given point $F$ in the plane and the usual polar coordinates $(\rho,\varphi)$ with the center $F$ introduce the inverse polar coordinates with the center $F$: $(r=\dfrac 1\rho$, $\varphi)$. We call the map $T$

\begin{equation}
T: \begin{cases}
\varphi_1=\varphi+ hr\\
r_1=r
\end{cases}
\end{equation}
the skew rotation with the center $F$ and with angular speed $h$.
\end{definition}


\begin{definition}
The map $T_h$ is called the perturbed skew rotation if it is invertible, symplectic and has the asymptotics
\begin{equation}
\label{eq: 1-2-3}
\begin{cases}
\varphi_1=\varphi+h r+ O(r^{2})\\
r_1=r+O(r^{3})
\end{cases}
\end{equation}
as $r\to 0$.
\end{definition}

Note that in this definition $(\varphi, r)\in [0,2\pi)\times \mathbb{R}_+$ are not necessarily the inverse polar coordinates but any "polar"  coordinate system in the neighbourhood of the infinity point $r=0$. If the correction terms $O(r^2)$, $O(r^3)$ on the right-hand side of \eqref{eq: 1-2-3} are absent, we call the map $T_h$ a generalized skew rotation.

%
\noindent\textbf{Remark.}
As it is known (see for example \cite{Arn}) a continuous time Hamiltonian system in the plane with closed trajectories can be written in the canonical form 
\begin{equation}
\label{eq: canonical}
\left\{
\begin{array}{cl}
\dot{\varphi}&=\omega(I)\\
\dot{I}&=0
\end{array}
\right.
\end{equation}
for particular variables $(\varphi, I)\in [0,2\pi)\times \mathbb{R}_+$ which are called action-angle coordinates. 
 The function $\omega(I)$ is called frequency or angular velocity. The action variable $I$ represents, up to the numerical factor, the area surrounded by the level curves of Hamiltonian.  Let the frequency  $\omega(I)$ monotonously tend to $0$ as $I\to \infty$. So we can use $(\omega,\varphi)$ as "polar" coordinates in a neighbourhood of infinity. In this coordinates, a shift along the orbits of the flow \eqref{eq: canonical} is the generalized skew rotation.

Let $T_1$ and $T_2$ be perturbed skew rotations in two different coordinate systems $(r,\varphi)$ and $(\tilde{r}, \tilde{\varphi})$ in a neighbourhood of infinity. What can one say about the existence of closed invariant curves for the map $T=T_1\cdot T_2$?

As we shall show, the qualitative behaviour of the product of several skew rotations is similar to the behaviour of the product of the usual rotations of the plane. Consider $N$ motions of the plane given, in complex notation, as $z\mapsto a_iz+b_i$, $|a_i|=1$. The product of these motions is \[z\mapsto az+b, \qquad a=\prod\limits_{i=1}^N a_i,\, b=\prod\limits_{i=1}^N b_i+\sum\limits_{i=0}^Nb_{N-i}\prod\limits_{j<i}a_{N-j}\] 
If $a\ne 1$ then this map is a rotation and so the trajectories of this map lie on the circles with the center $z=-\dfrac b a$ in the fixed point of this rotation. We generalize this result to the skew rotations and perturbed skew rotations. Our result shows that the typical trajectories lie on the closed invariant curves. The rest of the trajectories are contained in the annuli surrounded by these curves. This follows from the topological lemma.

\begin{lemma}
\label{lm: topological}
Let $A$ be the annulus bounded by two closed curves $\gamma_0$ and $\gamma_1$. Let $T$ be a homeomorphism of the domain $B\supset A$.
If $T(\gamma_i)=\gamma_i$, $i=0,1$, then $T(A)=A$.
\end{lemma}
  
On the heuristic level, our arguments are the following. A rotation around any center on the projective plane can be considered also as a rotation around infinite point. Thus the product of two rotations with different centres is a small perturbation of a rotation around infinity. For such a perturbed system with some non-degeneracy conditions, one can apply  Kolmogorov-Arnold-Moser theory and solve the described problem. 

However some work is needed to make this scheme rigorous.
We use a theorem  formulated in other notations in \cite[\S34]{Moser}.

\begin{theorem}[J. Moser]
\label{th: Moser}
Let $T$ be a real analytic invertible map defined in some "polar" coordinates $(\varphi,r)\in [0,2\pi)\times\mathbb{R}_+$ as 
\begin{equation}
\label{eq: Mosers_l+1}
\begin{cases}
\varphi_1=\varphi+h r^{2\ell}+ O(r^{2\ell+1}),\qquad h\ne 0\\
r_1=r+O(r^{2\ell+2})
\end{cases}
\end{equation}
for some $\ell\in \frac 12 \mathbb{N}$. Suppose also that $T$ satisfies the intersection property: the image of any closed curve surrounding the point $r=0$ and sufficiently close to the origin
intersects the original curve.

Then $T$ has an invariant curve surrounding the stable point $r=0$ in any neighbourhood of this point.

The measure $d\varphi dr$ of all points that do not belong to the closed invariant curves of $T$ is $o(\varepsilon)$ in the $\varepsilon$--neighbourhood of the stable point $r=0$.  
\end{theorem}

Theorem \ref{th: Moser} with $2\ell=1$ can be applied  to the perturbed skew rotations, due to the following lemma.

\begin{lemma}
\label{lm: symplectic}
Any perturbed skew rotation $T$ satisfies the intersection property.
\end{lemma}


From the Lemma \ref{lm: symplectic} and Theorem \ref{th: Moser}, it follows that the perturbed skew rotation has infinitely many closed invariant curves in any neighbourhood of infinity. By Lemma \ref{lm: topological} it means that the trajectory of any point is contained in some invariant annulus for the perturbed skew rotation.

In the next section we present the proof of the following theorem.

\begin{theorem}
\label{th: Rotations}
Let $R_{j,h_j}$ denote skew rotations with centres $F_j$, $j=1,\ldots,N$. Then for any set $(h_1,\ldots,h_N)$ such that $\sum\limits_{j=1}^N h_j\ne 0$, the trajectories of the map $T=\prod\limits_{j=1}^N R_{j,h_j}$ are bounded. Actually, they are contained in invariant annuli (as it is mentioned above). 
\end{theorem}

The case $\sum\limits_{j=1}^N h_j=0$ is considered in section \ref{sec: unbounded} (for $N=2$). In section \ref{sec: general} we present natural generalisations of Theorem \ref{th: Rotations}.  In section \ref{sec: squares} we analyze the non-smooth case. Finally, in section \ref{sec: numerics} we present the results of numerical computations.

Proofs of Lemmas \ref{lm: symplectic} and \ref{lm: topological} are given in Appendix. 

\paragraph*{Acknowledgements}
The authors thank Ya.G.~Sinai for many fruitful discussions and for his substantial help in writing of this article. We are also grateful to M. Blank for useful remarks.

\section{Proof of Theorem \ref{th: Rotations}.}
\label{sec: proof}

Since the product of two perturbed skew rotations in the same coordinate system is again a perturbed skew rotation, to prove Theorem \ref{th: Rotations} it is sufficient to show that the skew rotation with one center is the perturbed skew rotation in the inverse polar coordinate system with any other center. We verify the last statement by direct calculation. 

Without loss of generality consider the case of the skew rotation with the center $F_0=0$ in  the inverse polar coordinate system with the center $F_1=1$.  

The skew rotation $R_{0,h_0}$ is the shift along the trajectories of the system

\begin{equation}
\label{eq: skew_rotation_diff}
\begin{cases}
\dot{\varphi}=r\\
\dot{r}=0
\end{cases}
\end{equation}
where $(\varphi,r)$ are inverse polar coordinates with the center $F_0$. In the complex notation $\varphi=\mathrm{arg}(z)$, $r=\dfrac {1}{|z|}$. Inverse polar coordinates $(\tilde{\varphi},\tilde{r})$ with the center $F_1$ can be expressed as 

\begin{equation}
\begin{split}
\tilde{\varphi} = \mathrm{arg}(z-1)=\mathrm{Im}\left(\ln (z-1)\right) = \mathrm{Im} \left(\ln z+\ln \left(1-\frac1z\right)\right)= \\=\mathrm{Im} \left(\ln z -\frac1z +O(z^{-2})\right)=\varphi+r\sin \varphi +O(r^2)
\end{split}
\end{equation}

\begin{equation}
\label{eq: tilde_r}
\begin{split}
\tilde{r} =\frac{1}{|z-1|}=\exp\left(-\mathrm{Re}(\ln (z-1))\right) = \exp\left(-\mathrm{Re}( \ln z)- \mathrm{Re}\left(\ln \left(1-\frac1z\right)\right)\right)= \\=\frac{1}{|z|} \left(1 +\mathrm{Re}\left(\frac{1}{z}\right) +O(z^{-2})\right)=r(1+r\cos \varphi) +O(r^3)
\end{split}
\end{equation}

Thus system \eqref{eq: skew_rotation_diff} in coordinates $(\tilde{\varphi},\tilde{r})$ assumes the form

\begin{equation}
\label{eq: skew_rotation_tilde}
\begin{cases}
\dot{\tilde{\varphi}}=r(1+r\cos\varphi + O(r^2))=\tilde{r}+O(r^3)\\
\dot{\tilde{r}}=(-r^2\sin\varphi+O(r^3))r=O(r^3)
\end{cases}
\end{equation}

From \eqref{eq: tilde_r} it follows that $\dfrac{\tilde{r}}{r}\to 1$ as $r\to 0$, which gives $O(r^3)=O(\tilde{r}^3)$ in the neighbourhood of $r=0$. Since the map $R_{0,h_0}$ is a shift along the trajectories of \eqref{eq: skew_rotation_diff} it is also a perturbed skew rotation in the coordinate system $(\tilde{\varphi},\tilde{r})$. \ep

\section{The unbounded case.}
\label{sec: unbounded}
Here we analyze the case of products of two opposite skew rotations.

Recall the notion of the composition of two Hamiltonian systems with Hamiltonians $H_1$ and $H_2$ as a system corresponding to the vector field which is equal to a linear combination of vector fields given by $H_j$, $j=1,2$ (see \cite{Arn}). Such a system also can be written in Hamiltonian form with the Hamiltonian
\begin{equation}
\label{eq: Linear_combination}
H=h_1H_1+h_2H_2.
\end{equation} 
In the case of our basic example \eqref{eq: simple_shift} for $h_1\ne -h_2$, the level curves of the Hamiltonian $H$ are closed and are called Cartesian Ovals  -- the algebraic curves which are defined by  the equation 
\begin{equation}
\label{eq: Ovale}
h_1|z-F_1|+h_2|z-F_2|=const.
\end{equation}
For the particular case $h_1=h_2$, the level curves of the Hamiltonian $H$ is the family of confocal ellipses with foci $F_1$, $F_2$.
If $h_1=-h_2$, the relation \eqref{eq: Ovale} defines the family of confocal hyperbolas with the same foci. 

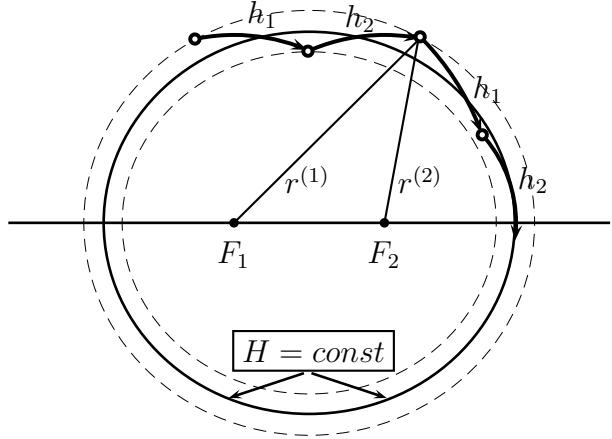
\begin{figure}[ht]
\label{fig: composition-mixture}
\begin{center}
\begin{pspicture}(0,0)(7,7)
\psline[linewidth=1 pt,linecolor=black,linestyle=solid]
(-0.5,3.5)(7.5,3.5)

\psdot(2.5,3.5)
\psdot(4.5,3.5)

\psellipse[linecolor=black,linewidth = 1 pt](3.5,3.5)(2.75,2.56)

\psellipse[linecolor=black,linewidth = 0.2 pt, linestyle=dashed](3.5,3.5)(2.49,2.28)

\psellipse[linecolor=black,linewidth = 0.2 pt, linestyle=dashed](3.5,3.5)(3,2.83)

\psarc[linewidth= 1.5 pt, linecolor=black]{<-o}(2.5,3.5){2.5}{67}{102}
\uput[12](2.5,6.2){$h_1$}

\psarc[linewidth= 1.5 pt, linecolor=black]{<-o}(4.5,3.5){2.5}{79}{114}
\uput[12](3.8,6.15){$h_2$}

\psarc[linewidth= 1.5 pt, linecolor=black]{<-o}(2.5,3.5){3.5}{20}{45}
\uput[12](5.5,5.2){$h_1$}

\psarc[linewidth= 1.5 pt, linecolor=black]{<-o}(4.5,3.5){1.75}{-8}{42}
\uput[12](6.1,4){$h_2$}

\psline[linecolor=black](2.5,3.5)(4.95,5.93)
\psline[linecolor=black](4.5,3.5)(4.95,5.93)

\uput[12](2.1,3){$F_1$}
\uput[12](4.1,3){$F_2$}

\uput[12](3,4){$r^{(1)}$}
\uput[12](4.5,4){$r^{(2)}$}
\uput[12](2.3,1.7){\psframebox[fillcolor=white]{$H=const$}}
\psline[linewidth= 1 pt, linecolor=black]{->}(3.43,1.5)(2.43,1.17)
\psline[linewidth= 1 pt, linecolor=black]{->}(3.53,1.5)(4.53,1.17)

\end{pspicture}
\end{center}
\caption{Linear combination of Hamiltonians $H_1$ and $H_2$.}
\end{figure}

In this context Theorem \ref{th: Rotations} states that for $h_1+h_2\ne 0$ the orbits of the map $T^{(h_1,h_2)}$ lie near the orbits of the completely integrable Hamiltonian system with Hamiltonian \eqref{eq: Linear_combination}.  For the case $h_1=-h_2$ the map $T^{(h,-h)}$ is an $O(|z|^{-3})$ perturbation of the identity map and thus KAM -- theory methods are not applicable. 

Let  $F_1=(-1,0)$ and $F_2=(1,0)$ be the two centres. Let $R_{1,h}$ and $R_{2,-h}$ be the skew rotations with centres $F_1$, $F_2$ with opposite values of parameters. One has the following result.

\begin{theorem}
\label{th: hyperbola}
The map $T=R_{1,h}\cdot R_{2,-h}$ has a trajectory escaping to infinity.
\end{theorem}

Informally, Theorem \ref{th: hyperbola} states that, in the case $h_1=-h_2$, the trajectories of the map $T^{(h_1,h_2)}$ also follow trajectories of the Hamiltonian system with Hamiltonian \eqref{eq: Linear_combination}.

\noindent\textbf{Proof.} 
Consider an arc $z_0\frown z_1$ of length $h$ with center $F_1$ starting at some point $z_0=x_0+i y_0$, $x_0>0$, which intersects the imaginary axis in its middle point $z_0^*=iy_0^*$. The distance from $z_1$ to the imaginary axis equals $-\mathrm{Re}(z_1)=-x_1=x_0$. The distance from  $z_2=R_2(z_1)$ to the axis is greater than $x_0$ since $\mathrm{Re}(z_1)<0$ and so $|z_1-F_2|>|z_0-F_2|$. Thus $\mathrm{Re}(z_2)>\mathrm{Re}(z_0)$. Similarly, for $z_3=R_1(z_2)$, we get $\mathrm{Re}(z_3)<\mathrm{Re}(z_1)$. 

By the same argument the distance from the $n$-th image of the point $z_0$ to the imaginary axis is greater then that of the $n-1$-st step. Incidentally, the arcs obtained on each iteration intersect the axis since $|x_j|<h$. Thanks to convexity of the circle, the sequence $\{y_j^*\}$ of coordinates of the intersection of the arc $z_j\frown z_{j+1}$ with the line $\mathrm{Re}(z)=0$ is increasing. Suppose, that this sequence has some finite limiting point $y$. Then for $(0,y)$ one gets $T(0,y)=(0, y)$ which is definitely impossible for any finite $y$. 

Similar considerations apply to any point $z\in \mathbb{C}$. The role of the axis $x=0$ in this case will be played by some hyperbola with foci $F_j$.  Theorem \ref{th: hyperbola} is proven.  \ep

\vspace{1 cm}
\begin{figure}[ht]
\begin{center}
\includegraphics[scale=0.4]{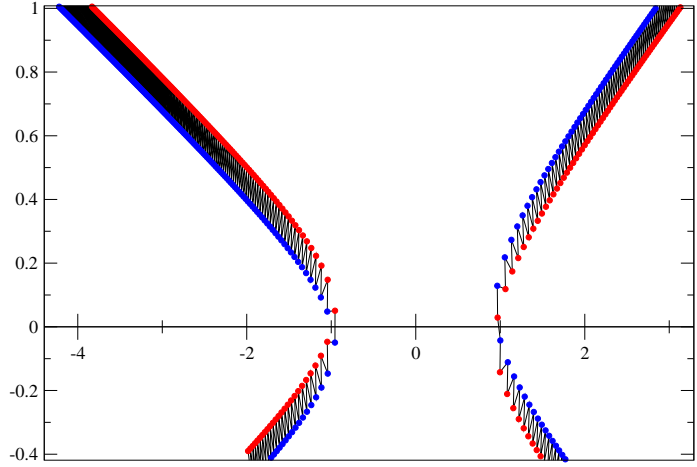} 
\end{center}
\caption{Hyperbolic trajectories.}
\end{figure}

 \section{Generalisations.}
 \label{sec: general}
 Now let us turn to the general case of a perturbed skew rotations given in a different coordinate systems. 
 \begin{definition}
Two coordinate systems $(r,\varphi)$ and $(\tilde{r},\tilde{\varphi})$ are called concordant in the neighbourhood of  $r=0$ if
\begin{enumerate}
\item The point $r=0$, $\varphi\in [0,2\pi)$ corresponds to $\tilde{r}=0$. So $\left.\dfrac{\partial \tilde{r}}{\partial \varphi}\right|_{r=0}=0$.
\item The transformation $(r,\varphi)\to(\tilde{r},\tilde{\varphi})$ is a real analytic invertible map.
\item This transformation asymptotically does not change the angles and radii, i.e.,  $\left.\dfrac{\partial \tilde{\varphi}}{\partial \varphi}\right|_{r=0}=1$ and \hspace {0.25 cm} $\left.\dfrac{\partial \tilde{r}}{\partial r}\right|_{r=0}=1$
\end{enumerate}
\end{definition}
Obviously the concordance of coordinate systems is an equivalence relation.

\noindent
\textbf{Example.}
Inverse polar coordinate systems with different centres are concordant.

\begin{lemma}
\label{lm: equithety}
If coordinate systems $(r,\varphi)$ and $(\tilde{r},\tilde{\varphi})$ are concordant in a neighbourhood of $r=0$ then $\left.\dfrac{\partial \tilde{r}}{\partial \varphi}\right|_{r=0}=O(r^2)$ for $r\to 0$.
\end{lemma}

\noindent\textbf{Proof.} From the definition $\left.\dfrac{\partial \tilde{r}}{\partial \varphi}\right|_{r=0}=0$ and $\left.\dfrac{\partial}{\partial r}\left(\dfrac{\partial \tilde{r}}{\partial \varphi}\right)\right|_{r=0}=\left.\dfrac{\partial}{\partial \varphi}\left(\dfrac{\partial \tilde{r}}{\partial r}\right)\right|_{r=0}=0$. \ep

\begin{theorem}
\label{th: equithety_inv}
Let $(r,\varphi)$ and $(\tilde{r},\tilde{\varphi})$ be two coordinate systems concordant in a neighbourhood of $r=0$. Let a real analytic invertible map $T$ satisfy conditions \eqref{eq: 1-2-3} in coordinates $(r,\varphi)$. Then $T$ also satisfies conditions \eqref{eq: 1-2-3} in coordinates $(\tilde{r},\tilde{\varphi})$.
\end{theorem}

\noindent\textbf{Proof.}
First consider 
\[\tilde{r}_1-\tilde{r}=\dfrac{\partial \tilde{r}}{\partial r} (r_1-r)+\dfrac{\partial \tilde{r}}{\partial \varphi} (\varphi_1-\varphi).\]
In a neighbourhood of $r=0$, the function $\dfrac{\partial \tilde{r}}{\partial r}$ is bounded and the factor $(r_1-r)$ is $O(r^3)$, due to conditions \eqref{eq: 1-2-3}.  For the second term, we get $\dfrac{\partial \tilde{r}}{\partial \varphi}=O(r^2)$ from 
Lemma \ref{lm: equithety} and $(\varphi_1-\varphi)=O(r)$, due to  conditions \eqref{eq: 1-2-3}. So we have 
\[\tilde{r}_1-\tilde{r}=O(r^3)=O(\tilde{r}^3)\]
since $\tilde{r}\sim r$.

For the difference of angular coordinates we have
\[\tilde{\varphi}_1-\tilde{\varphi}=\dfrac{\partial \tilde{\varphi}}{\partial r} (r_1-r)+\dfrac{\partial \tilde{\varphi}}{\partial \varphi} (\varphi_1-\varphi).\]
The first term is $O(r^3)$ from the previous discussion, and the second term is 
\[\dfrac{\partial \tilde{\varphi}}{\partial \varphi} (\varphi_1-\varphi)=(1+O(r))(h r+O(r^2))=hr +O(r^2)=h\tilde{r} +O(\tilde{r}^2).\]
Theorem \ref{th: equithety_inv} is proven. \ep

Now we can generalize Theorem \ref{th: Rotations}.

\begin{theorem}
\label{th: Hurray}
Let $(r,\varphi)$ and $(\tilde{r},\tilde{\varphi})$ be two coordinate systems concordant in a neighbourhood of $r=0$. Let a real analytic invertible map $T_1$  be a perturbed skew rotation in coordinates $(r,\varphi)$, a real analytic invertible map $T_2$  be a perturbed skew rotation in coordinates $(\tilde{r},\tilde{\varphi})$.  Then the product $T=T_1T_2$ is a perturbed skew rotation in coordinates $(r,\varphi)$ (or $(\tilde{r},\tilde{\varphi})$ as well). If the sum of the "angular coefficients" does not equal to zero then the conclusion of Theorem \ref{th: Moser} holds for $T$.
\end{theorem} 

\noindent\textbf{Proof.} It follows from Lemma \ref{lm: symplectic}, Theorem \ref{th: equithety_inv} and Theorem \ref{th: Moser}.\ep

\begin{corollary}
If there are two Hamiltonian systems with closed trajectories such that the functions $\omega_1$, $\omega_2$ monotonously tend to $0$ at infinity, the coordinate systems $(\varphi_1,\omega_1)$ and $(\varphi_2,\omega_2)$ are concordant and $h_1+h_2\ne 0$ then for the product of $h_1$-- and $h_2$--shifts along the trajectories of these Hamiltonian systems the conclusion of Theorem \ref{th: Moser}  holds.
\end{corollary}

\noindent\textbf{Remark.}
Theorem \ref{th: Hurray}  can be easily generalized to the case of any finite collection of invertible maps $T_1,\ldots, T_N$.

\section{Non-smooth case.}
\label{sec: squares}

In this section we show that the smoothness of the above discussed system is essential. Now we orient the plane opposite to the usual way. So the positive rotation direction is clockwise. Let $O_1=(-\frac 12,0)$ and $O_2=(\frac 12,0)$. Recall the basic example of the product of two skew rotations with different centres. Transformation $T$ moves each point $(x,y)$ along the level curve of the first Hamiltonian for the distance $h$ and then moves its image along the level curve of the second Hamiltonian for the distance $h$. In the basic example such level curves form two one-parametric families of concentric circles with centres $O_1$ and $O_2$. Consider the analogous system with two families of concentric squares, i.e. $H_1=|x-x_1|+|y-y_1|$ and $H_2=|x-x_2|+|y-y_2|$. For this system we show the coexistence of countably many periodic orbits and countably many trajectories escaping to infinity with speed $O(\sqrt{t})$.

\begin{figure}[ht]
\begin{center}
\begin{pspicture}(0,0)(5,5)

\psline[linewidth=1 pt,linecolor=black,linestyle=solid]
(1,0)(1,5)
\psline[linewidth=1 pt,linecolor=black,linestyle=solid]
(4,0)(4,5)

\psline[linewidth=0.5 pt,linecolor=black,linestyle=dotted, arrows=>->]
(-1.2,2.5)(0,3.7)
\psline[linewidth=0.5 pt,linecolor=black,linestyle=dotted, arrows=->]
(0,3.7)(1,4.7)

\psline[linewidth=0.5 pt,linecolor=black,linestyle=dotted, arrows=>->]
(-0.8,2.5)(0.2,3.5)
\psline[linewidth=0.5 pt,linecolor=black,linestyle=dotted, arrows=->]
(0.2,3.5)(1,4.3)

\psline[linewidth=0.5 pt,linecolor=black,linestyle=dotted, arrows=>->]
(-0.4,2.5)(0.4,3.3)
\psline[linewidth=0.5 pt,linecolor=black,linestyle=dotted, arrows=->]
(0.4,3.3)(1,3.9)

\psline[linewidth=0.5 pt,linecolor=black,linestyle=dotted, arrows=->]
(0,1.3)(-1.2,2.5)
\psline[linewidth=0.5 pt,linecolor=black,linestyle=dotted, arrows=>->]
(1,0.3)(0,1.3)

\psline[linewidth=0.5 pt,linecolor=black,linestyle=dotted, arrows=->]
(0,1.7)(-0.8,2.5)
\psline[linewidth=0.5 pt,linecolor=black,linestyle=dotted, arrows=>->]
(1,0.7)(0,1.7)

\psline[linewidth=0.5 pt,linecolor=black,linestyle=dotted, arrows=->]
(0,2.1)(-0.4,2.5)
\psline[linewidth=0.5 pt,linecolor=black,linestyle=dotted, arrows=>->]
(1,1.1)(0,2.1)

\psline[linewidth=0.5 pt,linecolor=black,linestyle=dotted, arrows=<-<]
(6.2,2.5)(5,3.7)
\psline[linewidth=0.5 pt,linecolor=black,linestyle=dotted, arrows=-<]
(5,3.7)(4,4.7)

\psline[linewidth=0.5 pt,linecolor=black,linestyle=dotted, arrows=<-<]
(5.8,2.5)(4.8,3.5)
\psline[linewidth=0.5 pt,linecolor=black,linestyle=dotted, arrows=-<]
(4.8,3.5)(4,4.3)

\psline[linewidth=0.5 pt,linecolor=black,linestyle=dotted, arrows=<-<]
(5.4,2.5)(4.6,3.3)
\psline[linewidth=0.5 pt,linecolor=black,linestyle=dotted, arrows=-<]
(4.6,3.3)(4,3.9)

\psline[linewidth=0.5 pt,linecolor=black,linestyle=dotted, arrows=>->]
(6.2,2.5)(5,1.3)
\psline[linewidth=0.5 pt,linecolor=black,linestyle=dotted, arrows=->]
(5,1.3)(4,0.3)

\psline[linewidth=0.5 pt,linecolor=black,linestyle=dotted, arrows=>->]
(5.8,2.5)(4.8,1.5)
\psline[linewidth=0.5 pt,linecolor=black,linestyle=dotted, arrows=->]
(4.8,1.5)(4,0.7)

\psline[linewidth=0.5 pt,linecolor=black,linestyle=dotted, arrows=>->]
(5.4,2.5)(4.6,1.7)
\psline[linewidth=0.5 pt,linecolor=black,linestyle=dotted, arrows=->]
(4.6,1.7)(4,1.1)
\pscustom[fillstyle=solid,fillcolor=gray!40,linestyle=none]
{
\psline(1,5)(1,2.5)
\psline(1,2.5)(4,2.5)
\psline(4,2.5)(4,5)
}
\pscustom[fillstyle=solid,fillcolor=gray!40,linestyle=none]
{
\psline(1,0)(1,2.5)
\psline(1,2.5)(4,2.5)
\psline(4,2.5)(4,0)
}
\psline[linewidth=1 pt,linecolor=black,linestyle=solid]
(-2,2.5)(7,2.5)
\uput[12](1,2.1){$O_1$}
\uput[12](3.3,2.1){$O_2$}

\psline[linewidth=0.5 pt,linecolor=red,linestyle=dashed, arrows=>->]
(1,3)(2,4)
\psline[linewidth=0.5 pt,linecolor=blue,linestyle=dashed, arrows=>->]
(2,4)(3,3)
\psline[linewidth=0.5 pt,linecolor=red,linestyle=dashed, arrows=>->]
(3,3)(4,4)

\psline[linewidth=0.5 pt,linecolor=blue,linestyle=dashed, arrows=>->]
(1,4)(2,3)

\psline[linewidth=0.5 pt,linecolor=red,linestyle=dashed, arrows=>->]
(2,3)(3,4)

\psline[linewidth=0.5 pt,linecolor=blue,linestyle=dashed, arrows=>->]
(3,4)(4,3)

\end{pspicture}
\end{center}
\label{fig: regions}
\caption{Three regions.}
\end{figure}
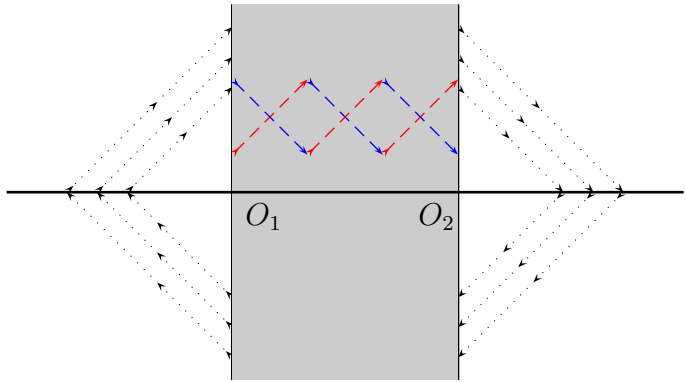

Let $Q_R^{(j)}=\{(x,y): |x-x_j|+|y-y_j|=R\}$, $j=1,2$ be two families of concentric squares with centres $O_1$ and $O_2$ such that the corresponding sides of the squares are parallel and $O_1O_2$ belongs to the common diagonal of the squares. The trajectory consists of odd and even steps. For even steps the point is shifted along the side of the square of the first family for the distance $h$ and for odd steps -- along the side of the square from the second family for the distance $h$ in the same direction. The plane is splitted onto three parts (see Fig. \ref{fig: regions}): half-planes $x>\dfrac 12$, $x<-\dfrac 12$ and the strip $|x|<\dfrac 12$. In the half-planes both dynamics coincide, while in the strip there exist a non-zero angle between the curves of first and second family. 

Let $h=a\sqrt{2}$. The shift along the square $Q_{t}^{(j)}$  for the distance $h$ sufficiently far from the vertices of the square corresponds to the shift for the distance $a$ in each coordinate in the positive or negative direction, depending on the part of the square.  

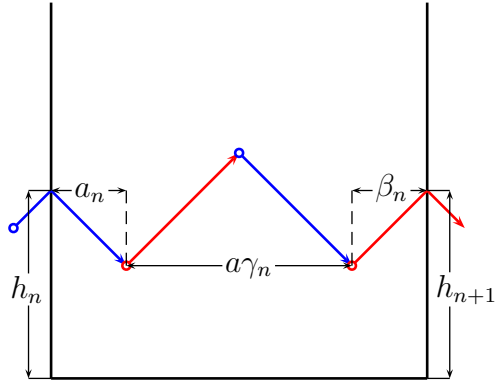
\begin{figure}[ht]
\begin{center}
\begin{pspicture}(0,0)(5,5)
\psline[linewidth=1 pt,linecolor=black,linestyle=solid]
(0,0)(0,5)

\psline[linewidth=1 pt,linecolor=black,linestyle=solid]
(0,0)(5,0)

\psline[linewidth=1 pt,linecolor=black,linestyle=solid]
(5,0)(5,5)

\psline[linewidth=1 pt,linecolor=blue,linestyle=solid, arrows = o-]
(-0.5,2)(0,2.5)
\psline[linewidth=1 pt,linecolor=blue,linestyle=solid, arrows = ->]
(0,2.5)(1,1.5)
\psline[linewidth=1 pt,linecolor=red,linestyle=solid, arrows = o->]
(1,1.5)(2.5,3)
\psline[linewidth=1 pt,linecolor=blue,linestyle=solid, arrows = o->]
(2.5,3)(4,1.5)
\psline[linewidth=1 pt,linecolor=red,linestyle=solid, arrows = o-]
(4,1.5)(5,2.5)
\psline[linewidth=1 pt,linecolor=red,linestyle=solid, arrows = ->]
(5,2.5)(5.5,2)

\psline[linewidth=0.5 pt,linecolor=black,linestyle=dashed]
(1,1.5)(1,2.5)

\psline[linewidth=0.5 pt,linecolor=black,linestyle=dashed]
(4,1.5)(4,2.5)

\psline[linewidth=0.5 pt,linecolor=black,linestyle=dashed, arrows = <->]
(4,2.5)(5,2.5)

\psline[linewidth=0.5 pt,linecolor=black,linestyle=dashed, arrows = <->]
(1,2.5)(0,2.5)

\psline[linewidth=0.5 pt,linecolor=black,linestyle=solid, arrows = <->]
(1,1.5)(4,1.5)
\psline[linewidth=0.5 pt,linecolor=black,linestyle=solid, arrows = <->]
(5.3,2.5)(5.3,0)
\psline[linewidth=0.5 pt,linecolor=black,linestyle=solid, arrows = <->]
(-0.3,2.5)(-0.3,0)
\psline[linewidth=0.5 pt,linecolor=black,linestyle=solid]
(5,2.5)(5.3,2.5)
\psline[linewidth=0.5 pt,linecolor=black,linestyle=solid]
(5,0)(5.3,0)
\psline[linewidth=0.5 pt,linecolor=black,linestyle=solid]
(0,2.5)(-0.3,2.5)
\psline[linewidth=0.5 pt,linecolor=black,linestyle=solid]
(0,0)(-0.3,0)
\psset{framesep=1 pt}
\uput[12](-0.75,1.1){\psframebox*{$h_n$}}
\uput[12](4.9,1.1){\psframebox*{$h_{n+1}$}}
\uput[12](0.1,2.4){\psframebox*{$a_n$}}
\uput[12](2.1,1.4){\psframebox*{$a\gamma_n$}}
\uput[12](4.1,2.45){\psframebox*{$\beta_n$}}
\end{pspicture}
\end{center}
\label{fig: shaft}
\caption{Dynamics in the strip.}
\end{figure}

Consider the trajectory of the point $(-\frac 12,h_0)$. Denote by $h_n$ the ordinate of the point of $n$-th entry to the strip, by $a_n\sqrt{2}$ the length of the remaining part of the shift inside the strip and by $\alpha_n$ the number of the family modulo $2$ (see Fig. \ref{fig: shaft}). 
For the starting point the corresponding triple is $(h_0,a,0)$. Whenever $|h_{n}|>a$, the values $(h_{n+1},a_{n+1},\alpha_{n+1})$ can be obtained by the following recurrent relations.

Denote by $\gamma_n$  the number of complete steps that the trajectory spends in the strip after the $n$-th entry and by $\beta_n\sqrt{2}$ the initial part of the next step belonging to the strip. We have
\begin{equation}
\label{eq: gamma_n}
\gamma_n=\left[\frac{1-a_n}{a}\right],
\end{equation}
\begin{equation}
\label{eq: beta_n}
\beta_n=1-a_n-a\gamma_n=a\left(\frac{1-a_n}{a}-\left[\frac{1-a_n}{a}\right]\right)=a\left\{\frac{1-a_n}{a}\right\},
\end{equation}
where $[x]$ and $\{x\}$ denote the integer and the fractional parts of $x$. For $h_{n+1}$ it immediately follows that
\begin{equation}
\label{eq: h_n}
h_{n+1}=-\left(h_n+(-1)^{\alpha_n}\left(a_n+(-1)^{\gamma_n+1}\beta_n-\frac{(1-(-1)^{\gamma_n})a}{2}\right)\right).
\end{equation}    

Compute the number of complete steps of the trajectory before the next entry to the strip. It is equal to $\left[\dfrac{2|h_{n+1}|+\beta_n-a}{a}\right]$.  Thus we get 
\begin{equation}
\label{eq: alpha_n}
\alpha_{n+1}\equiv \left(\alpha_n+\gamma_n+1+\left[\frac{2|h_{n+1}|+\beta_n}{a}\right]\right) \mathrm{mod}\, 2
\end{equation}
and 
\begin{equation}
\label{eq: a_n}
a_{n+1}=a\left(1-\left\{\frac{2|h_{n+1}|+\beta_n}{a}\right\}\right).
\end{equation}

In order to simplify the analysis let us consider the case $a=\dfrac{1}{m}$, $m\in \mathbb{N}$. 

\begin{proposition}
If $h_0=\ell a$, $\ell\in\frac 12\mathbb{N}$, $\ell\geqslant 1$ and $a_0=a$ then  $h_n=\ell_n a$, $\ell_n\in \frac 12\mathbb{N}$,  $\ell_n\geqslant 1$ and $a_n=a$ for all $n\in \mathbb{N}$.
\end{proposition}

\noindent\textbf{Proof} proceeds by induction. From $a_{n-1}=a$ and \eqref{eq: beta_n} it follows that $\beta_{n-1}=0$. Thus $h_{n}$ in \eqref{eq: h_n} is equal to $-h_{n-1}$ or $-h_{n-1}\pm a$ and so it is again of the form $\ell_n a$ where $2\ell_n \in \mathbb{N}$. Since $\left\{\dfrac{2|h_n|}{a}\right\}=0$ from \eqref{eq: a_n}, it follows that $a_n=a$. \ep 

Let $a=\dfrac {1}{2m}$, $a_0=a$ and $h_0$ be as above. Then  it follows from \eqref{eq: gamma_n} that $\gamma_n=2m-1$ and $h_{n+1}=-h_n$. So after $4(\ell+m)$ steps, the point $(-\frac 12, h_0)$ returns to the initial position.

For $a=\dfrac{1}{2m-1}$, $a_0=a$ it follows that $\gamma_n=2m-2$ and so $h_{n+1}=-h_n- (-1)^{\alpha_n}a$. 
So if $h_0=\ell a$, $\ell\in \mathbb{N}$ then $h_{2n+2}=h_{2n}+2a$ and the trajectory is expanding. If $h_0=\ell a$, $\ell\in \frac 12 \mathbb{N}\setminus \mathbb{N}$ then $h_{2n+2}=h_{2n}$ and the trajectory is periodic with the period $4\ell+2m-1$. More generally, if $a=\dfrac 1k$, $k\in \mathbb{N}$, $a_0=a-\varepsilon$,  $\varepsilon\in [0,a)$ then for $h_0=\ell a+\varepsilon$, $\ell\in \frac 12 \mathbb{N}\setminus \mathbb{N}$ the trajectory is again periodic.
This complete the proof of coexistence of periodic and escaping orbits.

Numerical experiments shows that for initial triples $(h_0,a,\alpha)$ with generic $h_0$ and $a$ the behaviour of the trajectory resembles a random walk (see Fig. \ref{fig: squares_gen}). 

\section{Numerical results.}
\label{sec: numerics}

First we present the results showing that for the case of two rotations the arguments from the section \ref{sec: unbounded} actually apply.
Computer simulations show that a trajectory of the map $T^{(h_1,h_2)}$ at each step intersects one level curve of the Hamiltonian \eqref{eq: Linear_combination}. 
\vspace{1 cm}
\begin{figure}[ht]
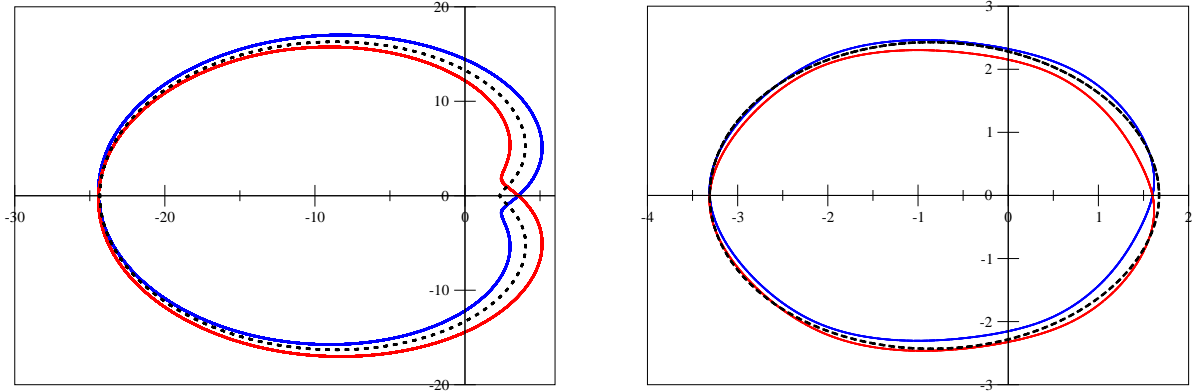

\includegraphics[width=0.46\textwidth]{oval4}\hspace{0.06\textwidth}
\includegraphics[width=0.46\textwidth]{oval_plus}
\caption{Left: Trajectory of the point $z=3+5i$ under the map $T^{(h_1,h_2)}$ for $h_1=2.5$, $h_2=-3$ intersects Cartesian Oval  $5|z-1|-6|z+1|=const$ (black dotted line depicts Cartesian Oval, upper grey curve consists of the points $T^n(S_1^{h_1} z)$, lower curve corresponds to $T^n(S_2^{h_2} z)$). Right:
Trajectory of the point $z=1.6+0i$ under the map $T^{(h_1,h_2)}$ for $h_1=2.5$, $h_2=0.25$ intersects Cartesian Oval $|z-1|+10|z+1|=const$. }
\label{fig: descartes}
\end{figure}

The map $T^{(h_1,h_2)}$ shows a typical KAM-theory behaviour: the islands of stability surrounded by closed invariant curves. 

\vspace{1 cm}
\begin{figure}[ht]
\centering
\includegraphics[width=0.5\textwidth]{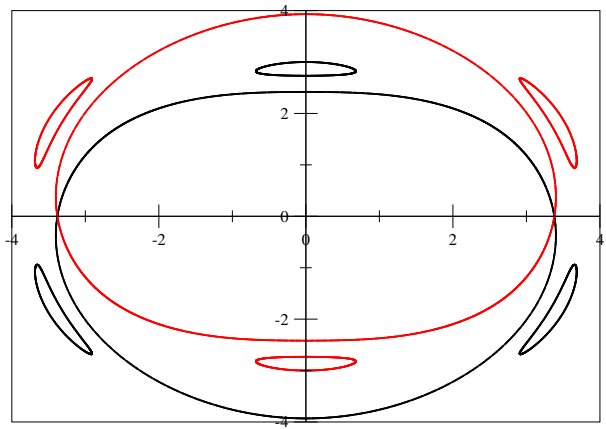}
\caption{Trajectories of the map $T^{(h_1,h_2)}$ for $h_1=h_2=3.8$ and initial points $(0,2.419)$ (large ovals) and $(0,3)$ (islands surrounding the periodic points).}
\end{figure}

\begin{figure}[ht]
\centering
\includegraphics[width=0.45\textwidth]{squares1}\hspace{0.06\textwidth}
\includegraphics[width=0.45\textwidth]{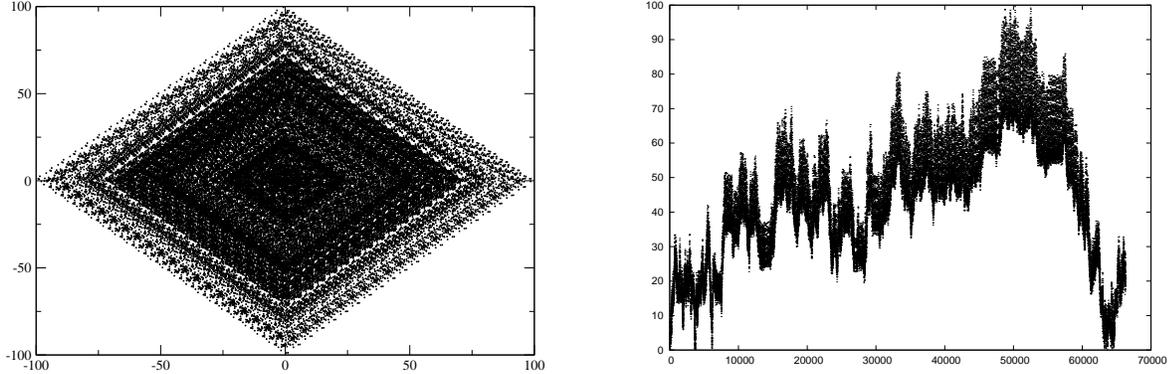}
\caption{Left: The points of the trajectory. Right: The distance of the trajectory $|T^n(x,y)|$ from $0$ for the step $h=2.43\sqrt{2}$ and randomly choosen initial point$(x_0,y_0)$: $|x_0|+|y_0|=1$.}
\label{fig: squares_gen}
\end{figure}

For the non-smooth case the dynamics looks like random walk. Typical trajectory is presented on figure  \ref{fig: squares_gen}.  
\vspace{1 cm}

\section*{Appendix.}
\noindent\textbf {Proof of Lemma \ref{lm: symplectic}.}  
The perturbed skew rotation $T$ is by definition an invertible area-preserving map. Consider any closed curve $\gamma$ surrounding points $z_1=0$ and $z_2=T^{-1}(0)$.  Then $T(\gamma)$ is a closed curve surrounding point $z_1=0$ and so, due to the area-preserving property, $T(\gamma)\cap \gamma\ne \emptyset$.  Since any closed curve surrounding the infinity in some of its neighbourhood also surrounds points $z_1$ and $z_2$, the intersection property  in the neighbourhood of infinity holds.\ep

\noindent\textbf{Proof of Lemma \ref{lm: topological}.} 
Suppose that for some $x_0\in A$ we have $T(x_0)\notin A$. Consider any curve $\gamma(t)\subset A$ having the following properties:
\begin{enumerate}
\item $\gamma(0)\in\gamma_0$, $\gamma(1)\in \gamma_1$,
\item \label{pr: contr} $\forall t\in (0,1):$ $\gamma(t)\cap (\gamma_1\cup \gamma_0)=\emptyset$,
\item $\exists t_0\in (0,1): $ $\gamma(t_0)=x_0$.
\end{enumerate}
Since $T$ is a homeomorphism and $T(\gamma(t_0))\notin A$, the image $T(\gamma(0,1))$ intersects the boundary $\gamma_1\cup \gamma_0$ in some point 
\[y_0=T(\gamma(t_1))\]
Then $\gamma(t_1)=T^{-1}(y_0)\in \gamma((0,1))\cap (\gamma_1\cup \gamma_0)$ which contradicts with our choice of $\gamma$ (property \ref{pr: contr}).\ep

 \end{document}